\documentstyle[12pt]{article}
\input{amssym.def}
\input{amssym.tex}
\begin{document}
\title{On complex and noncommutative tori}

\author{Igor  ~Nikolaev\\
Department of Mathematics\\
2500 University Drive N.W.\\
Calgary T2N 1N4 Canada\\
{\sf nikolaev@math.ucalgary.ca}}


 \maketitle


\newtheorem{thm}{Theorem}
\newtheorem{lem}{Lemma}
\newtheorem{dfn}{Definition}
\newtheorem{rmk}{Remark}
\newtheorem{cor}{Corollary}
\newtheorem{prp}{Proposition}
\newtheorem{exm}{Example}


\newcommand{\N}{{\Bbb N}}
\newcommand{\F}{{\cal F}}
\newcommand{\R}{{\Bbb R}}
\newcommand{\Z}{{\Bbb Z}}
\newcommand{\C}{{\Bbb C}}

\begin{abstract}
The ``noncommutative geometry'' of complex algebraic curves  is studied. 
As first step, we clarify a morphism between 
elliptic curves, or complex tori, and  $C^*$-algebras 
$T_{\theta}=\{u,v~|~vu=e^{2\pi i\theta}uv\}$, or noncommutative tori.
The main result says that under the morphism 
 isomorphic elliptic curves map to the Morita equivalent
noncommutative tori. Our approach is based on the rigidity of the length 
spectra of Riemann surfaces.

\vspace{7mm}


\vspace{5mm}
{\it AMS (MOS) Subj. Class.:  14H52, 46L85}
\end{abstract}

\section*{Introduction}
Noncommutative geometry is a branch of algebraic geometry studying 
``varieties'' over noncommutative rings. The noncommutative
rings are usually taken to be rings of operators acting on a Hilbert
space \cite{MvN}. The rudiments of noncommutative
geometry can be traced back to F.~Klein \cite{Kle1}, \cite{Kle2} or even earlier. 
The fundamental modern treatise \cite{C}  gives an account of 
status and perspective of the subject.

The noncommutative torus $T_{\theta}$ is a $C^*$-algebra generated by 
linear operators $u$ and $v$ on the Hilbert space $L^2(S^1)$
subjected to the commutation relation 
$vu=e^{2\pi i\theta}uv, ~\theta\in {\Bbb R}-{\Bbb Q}$  \cite{Rie}.
The classification of noncommutative tori was given in \cite{EfS}, \cite{PiV},
\cite{Rie}. Recall that two such tori $T_{\theta}, T_{\theta'}$ are Morita
equivalent if and only if $\theta,\theta'$ lie in the same orbit of the
action of group $GL(2,{\Bbb Z})$ on irrational numbers by linear fractional
transformations.

It is remarkable that ``moduli problem'' for $T_{\theta}$ looks
as such for the complex tori
$E_{\tau}={\Bbb C}/({\Bbb Z}+\tau {\Bbb Z})$, where $\tau$
is complex modulus. Namely, complex tori $E_{\tau},E_{\tau'}$ are isomorphic 
if and only if $\tau,\tau'$ lie in the same orbit of the action
of $SL(2,{\Bbb Z})$ on complex numbers by linear fractional
transformations. It was observed by some authors, e.g. \cite{Man}, \cite{Soi},
that it might be {\it not} just a coincidence.  This note is an attempt 
to show that it is indeed so: there exists  a general morphism  
between Riemann surfaces and $C^*$-algebras.

Let us give rough idea of our approach. Given Riemann surface $S$, there
is a function $S\to {\Bbb R}^{\infty}_+$ which maps the (discrete)
set of closed geodesics of $S$ to a discrete subset of real line
by assigning each closed geodesic its riemannian length. If $T_g(S)$
is the  space of all Riemann surfaces of genus $g\ge0$,
then function
\begin{equation}\label{eq1}
{\goth W}: T_g(S)\longrightarrow {\Bbb R}^{\infty}_+,
\end{equation}
is finitely-to-one and ``generically'' one-to-one \cite{Wol}. 
In the case $g=1$, function ${\goth W}$ is always
one-to-one. 
It is known also that restriction ${\goth W}_{syst}: T_g(S)\to {\Bbb R}_+$
of ${\goth W}$ to the shortest closed geodesic
of $S$ (called {\it systole})  is a $C^0$ Morse function on $T_g(S)$
\cite{SS}, \S 5.  Below we focus on the case $g=1$, i.e. $T_1\simeq E_{\tau}$.

Recall that $T_{\theta}$ has unique state $s_0$ (which is actually a tracial state)
\cite{Rie}. Any positive functional on $T_{\theta}$ has form $\omega s_0$,
where $\omega>0$ is a real number. Let $\Theta=\{T_{\theta}~|~\theta\in {\Bbb R}-{\Bbb Q}\}$
and $\Omega=\{\omega\in {\Bbb R}~|~\omega>0\}$. We define a map 
\begin{equation}\label{eq2}
{\goth V}:\Theta\times\Omega\longrightarrow {\Bbb R}^{\infty}_+
\end{equation}
by formula $(T_{\theta},\omega) \mapsto \{f_n(\omega)\ln tr (A_n)\}_{n=0}^{\infty}$,
where
\begin{eqnarray}
A_0 &=& \left(\small\matrix{a_0 & 1\cr 1 & 0}\right),\nonumber\\
A_1 &=& \left(\small\matrix{a_0 & 1\cr 1 & 0}\right) 
        \left(\small\matrix{a_1 & 1\cr 1 & 0}\right),\nonumber\\
    &\vdots& \\ 
A_n &=&  \left(\small\matrix{a_0 & 1\cr 1 & 0}\right) 
        \left(\small\matrix{a_1 & 1\cr 1 & 0}\right)\dots
        \left(\small\matrix{a_n & 1\cr 1 & 0}\right),\nonumber
\end{eqnarray}     	
are  integer  matrices whose entries $a_i>0$ are partial denominators
of continued fraction expansion of $\theta$ and $f_n$ are monotone $C^0$ functions
of $\omega$. Assuming that functions 
${\goth W}, {\goth V}$ have common range, one gets a mapping
${\goth W}{\goth V}^{-1}: E_{\tau}\to (T_{\theta},\omega)$.


Morphisms between $E_{\tau}$ and $T_{\theta}$ have
been studied in \cite{Man}, \cite{PoS}, \cite{Ram}, \cite{SoV}. 
The works \cite{Man}, \cite{Ram}
and \cite{SoV} treat noncommutative tori as 
``quantum compactification'' of the space of elliptic curves. 
This approach deals with an algebraic side 
of the subject. In particular,
Manin \cite{Man} suggested to use ``pseudolattices'' (i.e. $K_0$-group
of $T_{\theta}$) to solve the multiplication problem for real number    
fields. This problem  is part of the Hilbert 12th problem. In 
\cite{PoS} a functor from derived category of holomorphic vector bundles 
over $T_{\theta}$ to the Fukaya category of such bundles over $E_{\tau}$ 
was constructed. In this note we prove the  following results.  
\begin{thm}\label{thm1}
Let $E_{\tau}$ be complex torus of modulus $\tau, Im~\tau>0$
and $(T_{\theta},\omega)$ pair consisting of noncommutative 
torus with an irrational Rieffel's parameter $\theta$ and 
a positive functional $T_{\theta}\to {\Bbb C}$ of norm $\omega$. Then
there exists a one-to-one mapping $E_{\tau}\to (T_{\theta},\omega)$.
The action of modular group $SL(2,{\Bbb Z})$ on the complex
half plane $\{\tau\in {\Bbb C}~|~Im~\tau>0\}$ is equivariant
with:

\medskip
(i) the action of group $GL(2,{\Bbb Z})$ on irrationals
$\{\theta\in {\Bbb R}-{\Bbb Q}~|~\theta>0\}$ by linear
fractional transformations;

\smallskip
(ii) a discrete action on positive reals $\{\omega\in {\Bbb R}~|~\omega>0\}$.

\medskip\noindent
In particular, isomorphic complex tori map to the Morita equivalent
noncommutative tori, and vice versa.
\end{thm}
\begin{dfn}\label{dfn1}
The irrational number $\theta$ of  mapping $E_{\tau}\to (T_{\theta},\omega)$
we call a projective curvature of the elliptic curve $E_{\tau}$.
\end{dfn}
\begin{thm}\label{thm2}
Projective curvature of an elliptic curve with complex multiplication 
is a quadratic irrationality.
\end{thm}

\medskip\noindent
{\sf Acknowledgements.} It is my pleasure to thank G.~A.~Elliott, 
Yu.~I.~Manin   and  M.~Rieffel
for their interest in the subject of present note. I am greatful 
to the referee for critical remarks and helpful suggestions.   


\section{Proofs}
The proof of both theorems is based on the rigidity of length spectrum of complex
torus, cf Wolpert \cite{Wol}. A preliminary information on complex and
noncommutative tori can be found in Section 2.  
\subsection{Proof of Theorem 1}
Let us review main steps of the proof. By the rigidity lemma (Lemma \ref{lm1})
the length spectrum $Sp~ E$ defines conformal structure of $E$. 
In fact, this correspondence as a bijection. Under isomorphisms of $E$ the length
spectrum can acquire a real multiple or get a ``cut of finite tail''
(Lemma \ref{lm2}).
We attach to ${\Bbb C}/L$ a continued fraction of its projective curvature $\theta$
as specified in Introduction.
Then isomorphic 
tori ${\Bbb C}/L$ will have continued fractions which differ
only in a finite number of terms. In other words, one can attach a Morita equivalence 
class of noncommutative tori to every isomorphism class of complex tori.
\begin{lem}\label{lm1}
{\bf (Rigidity of length spectrum)}
Let $Sp~E$ be length spectrum of a complex torus $E={\Bbb C}/L$.
Then there exists a unique complex torus with the spectrum. This correspondence 
is a bijection.
\end{lem}
{\it Proof.}  See  McKean \cite{McK}.
$\square$

\bigskip\noindent
Let $Sp~X=\{l_1,l_2,\dots\}$ be length spectrum of a Riemann surface
$X$. Let $a>0$ be a real number. By $aSp~X$ we understand the length
spectrum $\{al_1,al_2,\dots\}$. Similarly, for any $m\in {\Bbb N}$ we
denote by $Sp_m~X$ the length spectrum $\{l_m,l_{m+1},\dots\}$, i.e. 
the one obtained by deleting the first $(m-1)$-geodesics in $Sp~X$. 
\begin{lem}\label{lm2}
Let $E\sim E'$ be isomorphic complex tori. Then either:

\medskip
(i) $Sp~E'=|\alpha|Sp~E$ for an $\alpha\in {\Bbb C}^{\times}$, or

\smallskip
(ii) $Sp~E'=Sp_m~E$ for a $m\in {\Bbb N}$.
\end{lem}
{\it Proof.}
(i) The complex tori $E={\Bbb C}/L, E'={\Bbb C}/M$ are
isomorphic if and only if $M=\alpha L$ for a complex number
$\alpha\in {\Bbb C}^{\times}$. It is not hard to see that 
closed geodesic of $E$ are bijective with the points of  the lattice 
$L=\omega_1{\Bbb Z}+\omega_2{\Bbb Z}$ in the following way. Take a
segment of straight line through points $0$ and $\omega$ of lattice
$L$ which contains no other points of $L$. This segment represents
a homotopy class of curves through $0$ and a closed geodesic of $E$. 
Evidently, this geodesic will be the shortest in its homotopy class with 
the length $|\omega|$ equal to absolute value of complex number $\omega$.    
Thus, $|\omega|$ belongs to the length spectrum of $E$.

Let now $Sp~E=\{|\omega_1|,|\omega_2|,\dots\}$ with $\omega_i\in L$. Since
$M=\alpha L$, one gets $Sp~E'=\{|\alpha||\omega_1|,|\alpha||\omega_2|,\dots\}$
and $Sp~E'=|\alpha|Sp~E$. Item (i) follows.

\medskip
(ii) Note that according to (i) the length spectrum $Sp~X=\{l_0,l_1,l_2,\dots\}$ can
be written as $Sp~X=\{1,l_1,l_2,\dots\}$ after multiplication on $1/l_0$, where $l_0$
is the length of shortest geodesic. Note also that shortest geodesic of 
complex torus has homotopy type $(1,0)$ or $(0,1)$ (standard generators for $\pi_1E$).

Let $a,b,c,d$ be integers such that $ad-bc=\pm 1$ and let 
\begin{eqnarray}\label{eq5}
\omega_1' &=& a\omega_1+b\omega_2,\nonumber\\
\omega_2' &=& c\omega_1+d\omega_2,
\end{eqnarray}
be an automorphism of the lattice $L=\omega_1{\Bbb Z}+\omega_2{\Bbb Z}$.
This automorphisms maps standard generators $(1,0)$ and $(0,1)$ of
$L$ to the vectors $\omega_1=(a,b),\omega_2=(c,d)$. Let their lengths
be $l_m,l_{m+1}$, respectively. 

As we showed earlier, $l_m,l_{m+1}\in Sp~E$ and it is not hard to see
that there are no geodesics of the intermediate length. (This gives a justification
for the notation chosen.) Note that $\omega_1,\omega_2$ are standard generators
for the complex torus $E'\sim E$ and therefore one of them is the shortest closed geodesics 
of $E'$. One can normalize it to the length $1$. 

On the other hand, there are only finite number of closed geodesics of length smaller than
$l_n$ (McKean \cite{McK}). Thus $Sp~E\cap Sp~E'=\{l_m,l_{m+1},\dots\}$ for a finite
number $m$ and since (\ref{eq5}) is automorphism of the lattice $L$. In other words, 
$Sp~E'=Sp_m~E$. Item (ii) follows.
$\square$

\medskip\noindent
To finish the proof of item (i) of Theorem \ref{thm1}, one needs to combine
Lemmas \ref{lm1},\ref{lm2} with the fact that two noncommutative
tori $T_{\theta},T_{\theta'}$ are Morita equivalent if and only
if their continued fractions differ only in a finite number
of terms (Section 2.1).

\medskip
To prove item (ii) of Theorem \ref{thm1}, let to the contrary the action
of $SL(2,{\Bbb Z})$ be non-discrete, i.e. having limit points in $\Omega$.
Let $p=\lim_{n\to\infty} (T_{\theta_n},\omega_n)$, where $\theta_n$ lie
in the same orbit of $GL(2,{\Bbb Z})$. Let $E_p$ be corresponding complex
torus such that $E_p\not\cong E_{p_n}$ are non-isomorphic. By continuity
of the systole function ${\goth W}_{syst}$ (see Introduction)
$Sp~E_p=\lim Sp ~E_{p_n}$. Then by the rigidity of length spectra,
$E_p\cong E_{p_n}$ are isomorphic. The contradiction proves item (ii)
of theorem. 
$\square$

\subsection{Proof of Theorem 2}
Let us outline the idea of the proof. If $E$ admits complex multiplication,
then its complex modulus $\tau$ lies in an imaginary quadratic field $K$.
In fact, up to an isogeny, the ring of endomorphisms $End~E={\cal O}_K$, where ${\cal O}_K$
is the ring of integers of field $K$. It can be shown that $L$ is an ideal
in ${\cal O}_K$ (Section 2.3). The length spectrum $Sp~E$ of elliptic curve
with complex multiplication is a ``geometric progression'' with the growth rate
$|\alpha|$, where $\alpha\in~End~E$ (Lemma \ref{lm5}). One can use Klein's
lemma (Lemma \ref{lm5b}) to characterize length spectra in terms of continued
fractions. In particular, length spectrum with asymptotically geometric 
growth correspond to periodic continued fractions. Thus, projective curvature 
converges to quadratic irrationality.
\begin{dfn}
Length spectrum $Sp~E$ of an elliptic curve $E={\Bbb C}/L$ is called
$\alpha$-multiplicative, if there exists a complex number $\alpha\in
{\Bbb C}^{\times}$ with $|\alpha|>1$ such that 
\begin{equation}
Sp~E=\{l_1,\dots,l_N,|\alpha|l_1,\dots,|\alpha|l_N,\dots,|\alpha|^nl_1,\dots,
|\alpha|^nl_N,\dots\},
\end{equation} 
\end{dfn}
for a $N\in{\Bbb N}$.
\begin{lem}\label{lm5}
Let $E$ be an elliptic curve with complex multiplication. Then 
its length spectrum $Sp~E$ is $\alpha$-multiplicative for an
$\alpha\in {\Bbb C}^{\times}$.
\end{lem}
{\it Proof.} Let $E={\Bbb C}/L$ be complex torus which admits non-trivial
endomorphisms $z\mapsto\alpha z,\alpha\in K={\Bbb Q}(\sqrt{-d})$. It is
known that $End~E$ is an order in the field $K$. In fact, up to an
isogeny of $E$, $End~E\simeq {\cal O}_K$, where ${\cal O}_K$ is the ring
of integers of imaginary quadratic field $K$ (Section 2.3). Lattice $L$
in this case corresponds to an ideal in ${\cal O}_K$.

Let $l_1$ be minimal length of closed geodesic of $E$. For an endomorphism
$\alpha: E\to E, \alpha\in {\Bbb C}^{\times}$, consider the set of 
geodesics whose lengths are less than $|\alpha|l_1$. By the properties
of $Sp~E$ mentioned in Section 1.1, such a set will be finite. Let us
denote the lengths of geodesics in this set by $l_1,\dots,l_N$. 
Since every geodesic in $Sp~E$ is a complex number $\omega_i$ 
lying in the ring $L\subseteq {\cal O}_K$, one can consider the set of geodesics 
$\alpha\omega_1,\dots,\alpha\omega_N$. The length of these geodesics will be 
$|\alpha|l_1,\dots,|\alpha|l_N$, respectively.  It is not hard to see that by 
the choice of number $N$, the first $2N$ elements of $Sp~E$ are presented by 
the following growing sequence of geodesics: $l_1,\dots,l_N, |\alpha|l_1,\dots,
|\alpha|l_N$.  We proceed by iterations of $\alpha$, until all closed geodesics of $E$ are
exhausted. The conclusion of  Lemma \ref{lm5} follows. 
$\square$

\medskip
We shall need the following statement regarding geometry of the regular 
continued fractions \cite{Kle1},\cite{Kle2}.
It is valid for any regular fraction, not necessarily periodic.
\begin{lem}\label{lm5b}
{\bf (F.~Klein)}
Let
\begin{equation}\label{eq9}
\omega=\mu_1+{1\over\displaystyle \mu_2+
{\strut 1\over\displaystyle \mu_3\displaystyle +\dots}}
\end{equation}
be a regular continued fraction. Let us denote  the convergents of $\omega$
by:
\begin{equation}\label{eq10} 
{p_{-1}\over q_{-1}}={0\over 1}, \quad
{p_{0}\over q_{0}}={1\over 0}, \quad
{p_{1}\over q_{1}}={\mu_1\over 1},\dots,
{p_{\nu}\over q_{\nu}}={\mu_{\nu}p_{\nu-1}+p_{\nu-2}\over
\mu_{\nu}q_{\nu-1}+q_{\nu-2}}.
\end{equation}
For any lattice $L$ in ${\Bbb C}$, consider a segment $I$ with ends in
the points  $(p_{\nu-2},q_{\nu-2})$ and  $(p_{\nu},q_{\nu})$. 
Then the segment $J$ which joins $0$ with the point $p_{\nu-1},q_{\nu-1}$
is parallel to $I$ and
\begin{equation}\label{eq10b}
|I|=\mu_{\nu}|J|,
\end{equation}
where $|\bullet|$ denotes the length of the segment.
\end{lem}
{\it Proof.} We refer the reader to \cite{Kle2}.
$\square$
\begin{cor}\label{cor1}
Let $\omega_{\nu}=(p_{\nu},q_{\nu})$ be lattice points
mentioned in Lemma \ref{lm5b}. Then the length 
of vector $\omega_{\nu}$ can be evaluated with the
help of the following asymptotic formula:
\begin{equation}\label{eq10c}
|\omega_{\nu}|\approx |\omega_{\nu-2}|+\mu_{\nu}|\omega_{\nu-1}|.
\end{equation}
\end{cor}
{\it Proof.} Indeed, using  notation of Lemma \ref{lm5b},
one can write $|(p_{\nu},q_{\nu})|\approx $\linebreak $ |(p_{\nu-2},q_{\nu-2})|
+|I|$. But according to equation (\ref{eq10b}), 
$|I|=\mu_{\nu}|(p_{\nu-1},q_{\nu-1})|$. Corollary
\ref{cor1} follows.
$\square$

\medskip
Note that according to the recurrent formula (\ref{eq10c}) the length spectrum
$\{|\omega_{\nu}|\}$ coming from continued fraction (\ref{eq9}) is completely
determined by the first two values: $|\omega_1|$ and $|\omega_2|$. Using (\ref{eq10c}),
one can easily deduce the following asymptotic formula for 
$|\omega_{\nu}|$ as function of $|\omega_1|,|\omega_2|$:
\begin{equation}\label{eq10d}
|\omega_{\nu}|\approx |\omega_2|\prod_{k=3}^{\nu}\mu_k+
|\omega_1|\prod_{k=4}^{\nu}\mu_k + O(\nu).
\end{equation}
Fix $N$ a positive integer. It follows from equation (\ref{eq10d})
that:
\begin{eqnarray}\label{eq10e} 
\lim_{\nu\to\infty}{|\omega_{\nu+N}|\over |\omega_{\nu}|} &=&
\mu_{\nu+1}\dots\mu_{\nu+N}
\lim_{\nu\to\infty}\left({\mu_{\nu}\dots\mu_3|\omega_2|+
\mu_{\nu}\dots\mu_4|\omega_1|+O(\nu)\over
\mu_{\nu}\dots\mu_3|\omega_2|+
\mu_{\nu}\dots\mu_4|\omega_1|+O(\nu)}\right)=\nonumber\\
 &=& \mu_{\nu+1}\dots\mu_{\nu+N}. 
\end{eqnarray}
Let $E$ be an elliptic curve with complex multiplication. Then
by Lemma \ref{lm5} its length spectrum $Sp~E$ is 
$\alpha$-multiplicative. In other words, 
\begin{equation}\label{eq10f}
{l_{\nu+N}\over l_{\nu}}=|\alpha|=~Const,
\end{equation}
for an $N\in {\Bbb N}$ and any $\nu~mod~N$. Note that
$|\alpha|$ is a rational integer. Thus, by formula
(\ref{eq10e}) we have $\mu_{\nu+1}\dots\mu_{\nu+N}=~Const$,
for any $\nu~mod~N$. The last requirement can be satisfied if 
and only if continued fraction (\ref{eq9}) is $N$-periodic.
Theorem \ref{thm2} is proven.
$\square$

\section{Background information}
In present section we briefly review noncommutative and complex tori.
The excellent source of information on noncommutative torus are papers
\cite{EfS}, \cite{Rie} and  monograph of  \cite{RLL}. The literature 
on complex torus is fairly vast. We recommend for the reference 
Ch. VI of \cite{S}.

\subsection{Noncommutative torus}
By the $C^*$-algebra one understands a noncommutative Banach
algebra with an involution \cite{RLL}. Namely, a $C^*$-algebra
$A$ is an algebra over $\Bbb C$ with a norm $a\mapsto ||a||$
and an involution $a\mapsto a^*, a\in A$, such that $A$ is
complete with respect to the norm, and such that 
$||ab||\le ||a||~||b||$ and $||a^*a||=||a||^2$ for every
$a,b\in A$. If $A$ is commutative, then the Gelfand
theorem says that $A$ is isometrically $*$-isomorphic   
to the $C^*$-algebra $C_0(X)$ of continuous complex-valued
functions on a locally compact Hausdorff space $X$.
For otherwise, $A$ represents a ``noncommutative'' topological
space $X$.

\medskip\noindent
\underline{$K_0$ and dimension groups.}
~Given a $C^*$-algebra, $A$, consider new $C^*$-algebra $M_n(A)$, 
i.e. the matrix
algebra over $A$. There exists a remarkable semi-group, $A^+$,
connected to the set of projections in algebra
$M_{\infty}=\cup_{n=1}^{\infty} M_n(A)$. Namely, projections 
$p,q\in M_{\infty}(A)$
are Murray-von Neumann equivalent $p\sim q$ if they can be presented as
$p=v^*v$ and $q=vv^*$ for an element $v\in M_{\infty}(A)$.
The equivalence class of projections is denoted by $[p]$.
The semi-group $A^+$ is defined to be the set of all equivalence
classes of projections in $M_{\infty}(A)$ with the binary operation
$[p]+[q]=[p\oplus q]$.
The Grothendieck completion of $A^+$ to an abelian group is called
a {\it $K_0$-group of $A$}.
The functor $A\to K_0(A)$ maps the unital
$C^*$-algebras into the category of abelian groups so that
the semi-group $A^+\subset A$ corresponds to a ``positive
cone'' $K_0^+\subset K_0(A)$ and the unit element
$1\in A$ corresponds to the ``order unit''
$[1]\in K_0(A)$.
The ordered abelian group $(K_0,K_0^+,[1])$ with the order
unit is called a {\it dimension (Elliott) group} of $A$. 
The dimension (Elliott) group is complete invariant of the $AF$
$C^*$-algebras.

\medskip\noindent
\underline{Noncommutative torus.}
~Fix $\theta$ irrational and consider a linear flow $\dot x=\theta,\dot y=1$
on the torus. Let $S^1$ be a closed transversal to our flow.
The {\it noncommutative torus} $T_{\theta}$
is a norm-closed $C^*$-algebra generated by the unitary operators
in the Hilbert space $L^2(S^1)$: 
\displaymath
Uf(t)=z(t)f(t),\qquad Vf(t)=f(t-\alpha),
\enddisplaymath
which are multiplication by a unimodular function $z(t)$
and rotation operators. It could be easily verified that $UV=e^{2\pi i\alpha}VU$.
As an ``abstract'' algebra, $T_{\theta}$ is a crossed product $C^*$-algebra 
$C(S^1)\rtimes_{\phi}{\Bbb Z}$ of (commutative) $C^*$-algebra of
complex-valued continuous functions on $S^1$ by the action of 
powers of $\phi$, where $\phi$ is a rotation of $S^1$ through
the angle $2\pi\alpha$. $T_{\theta}$ is not $AF$, but can be
embedded into an $AF$-algebra whose  dimension group is
$P_{\theta}$ (to be specified below);  the latter is known to be intimately
connected with the arithmetic of the irrational numbers $\theta$'s.
The following beautiful result is due to the efforts of many
mathematicians
\footnote{The author apologizes for possible erroneous credits 
regarding history of the problem. Classification of noncommutative
tori seems to be an old problem; early results in this direction 
can be found in the works of Klein \cite{Kle1},\cite{Kle2}.}     
(Effros, Elliott, Pimsner, Rieffel, Shen, Voiculescu, etc): 
\begin{thm} {\bf (Classification of noncommutative tori)}
Let $T_{\theta}$ be a noncommutative torus.
 Suppose that the $\theta$ has a continued fraction expansion 
\displaymath
\theta=a_0+{1\over\displaystyle a_1+
{\strut 1\over\displaystyle a_2\displaystyle +\dots}}
\buildrel\rm def\over= [a_0,a_1,a_2,\dots].
\enddisplaymath
Let $\varphi_n$ be a composition of isometries of the lattice
$\Z^2\subset\R^2:\quad\varphi_n=\left(
\small\matrix{ a_0 & 1\cr 1   & 0}\right)$\linebreak
$\dots\left(\small\matrix{ a_n & 1\cr
1   & 0}\right)$.
Then $T_{\theta}$ can be embedded into an $AF$-algebra whose
dimension group is a direct limit
of the ordered abelian groups:
$P_{\theta}=\lim_{n\to\infty}(\Z^2,\varphi_n)$.
Moreover, if $\theta=[a_0,a_1,\dots]$ and $\theta'=[b_0,b_1,\dots]$ are two
irrational numbers, then $P_{\theta}$ and $P_{\theta'}$
are isomorphic (i.e. noncommutative tori $T_{\theta}$ and $T_{\theta'}$ are
Morita equivalent) if and only if $a_{m+k}=b_m$ for an
integer number $k\in\Z$. In other words, the irrational numbers
$\theta$ and $\theta'$ are modular equivalent:
$\theta'={a\theta+b\over c\theta+d},\quad ad-bc=\pm 1$,
where $a,b,c,d\in\Z$ are integer numbers.
\end{thm}
{\it Proof.} An algebraic proof of this fact can be found in \cite{EfS}.
$\square$

\subsection{Complex torus}
Let $L$ denote a lattice in the complex plane ${\Bbb C}$. Attached to $L$,
there are following classic Weierstrass function $\wp(z;L)$ and Eisenstein 
series $G_k(L)$:
\begin{eqnarray}
\wp(z;L) &=& {1\over z^2}+\sum_{\omega\in L^{\times}}
\left\{{1\over (z+\omega)^2}-{1\over\omega^2}\right\},\\
G_k(L) &=& \sum_{\omega\in L^{\times}}^{k\ge2}{1\over\omega^{2k}}.
\end{eqnarray}
$\wp(z;L)$ is analytic and $G_k(L)$ is convergent for any lattice $L$
\cite{S}. There exists a duality between lattices $L$ and cubic curves $E$ 
given by the following theorem.
\begin{thm}\label{thm4}
Let $L$ be a lattice in $\Bbb C$. Then the map
$z\mapsto (\wp(z;L),{1\over 2}\wp'(z;L))$ is an analytic isomorphism
from complex torus ${\Bbb C}/L$ to elliptic cubic $E=E({\Bbb C})$:
\begin{equation}
E({\Bbb C})=\{(x,y)\in {\Bbb C}^2~|~y^2=x^3-15G_4(L)x-35G_6(L)\}.
\end{equation}
Conversely, to any cubic in the Weierstrass normal form
$y^2=x^3+ax+b$ there corresponds a unique lattice $L$ such
that $a=-15G_4(L)$ and $b=-35G_6(L)$.  
\end{thm}
{\it Proof.} We refer the reader to \cite{S} for a detailed
proof of this fact.
$\square$

\medskip\noindent
\underline{Isomorphism of complex tori.}
~Let $L$ be a lattice in $\Bbb C$. The Riemann surface
${\Bbb C}/L$ is called a {\it complex torus}. Let
$f:{\Bbb C}/L\to {\Bbb C}/M$ be holomorphic and invertible
map (isomorphism) between two complex tori. Since $f$
is covered by a linear map $z\to\alpha z$ on $\Bbb C$,
one can easily conclude that $\alpha L=M$ for an 
$\alpha\in {\Bbb C}^{\times}$. On the other hand,
lattice $L$ can always be written as $L=\omega_1{\Bbb Z}+\omega_2{\Bbb Z}$,
where $\omega_1,\omega_2\in {\Bbb C}^{\times}$ and 
$\omega_2\ne k\omega_1$ for a $k\in {\Bbb R}$. The complex
number $\tau={\omega_2\over\omega_1}$ is called a {\it complex
modulus} of lattice $L$.   
\begin{lem}\label{lm6}
Two complex tori are isomorphic if and only if their complex moduli
$\tau$ and $\tau'$ satisfy the equation:
\begin{equation}\label{eq17}
\tau'={a\tau+b\over c\tau+d}
\qquad ad-bc=\pm 1,\qquad a,b,c,d\in {\Bbb Z}.
\end{equation}
\end{lem}
{\it Proof.} The proof of this fact can be found in \cite{S}. 
$\square$

\subsection{Elliptic curves with complex multiplication}
Let $E={\Bbb C}/L$ be an elliptic curve. Consider the set
$End~E$ of analytic self-mappings of $E$. Each $f\in End~E$
is covered on the complex plane by map $z\mapsto\alpha z$ 
for an $\alpha\in {\Bbb C}$. It is not hard to see that $End~E$
has the structure of a ring under the pointwise addition and
multiplication of functions. The set $End~E$ is called
an {\it endomorphism ring} of elliptic curve $E$. By the remarks 
above, $End~E$ can be thought of as a subring of complex numbers:
\begin{equation}\label{eq20}
End~E=\{\alpha\in {\Bbb C}~|~\alpha L\subset L\}.
\end{equation}
There exists a fairly complete algebraic  description of such
rings. Roughly speaking, they are either 
``rational integers'' $\Bbb Z$ or integers ${\cal O}_K$ of an
algebraic number field $K$. The following lemma is true.   
\begin{lem}\label{lm8}
Let $\alpha\in End~E$ be a complex number. Then either:

\medskip
(i) $\alpha$ is a rational integer, or

\smallskip
(ii) $\alpha$ is an algebraic integer in an imaginary
quadratic number field $K={\Bbb Q}(\sqrt{-d})$. 
\end{lem}
{\it Proof.} See \cite{S}.
$\square$

\medskip\noindent
\underline{Complex multiplication.}
~If $End~E $ is different from $\Bbb Z$,
$E$ is said to be {\it elliptic curve with complex multiplication}. 
If $E$ admits complex multiplication, then its ring $End~E$ is
an order in an imaginary quadratic field $K$. In fact, $E$
admits an isogeny (analytic homomorphism) to a curve $E'$
such that $End~E'\simeq {\cal O}_K$, where ${\cal O}_K$
is the ring of integers of field $K$ \cite{S}. 
Thus, by property $\alpha L\subseteq L$, lattice
$L$ is an ideal in ${\cal O}_K$. Denote by $h_K$ the
class number of field $K$. It is well known, that
there exist $h_K$ non-isomorphic ideals in ${\cal O}_K$.
Therefore, elliptic curves $E_1={\Bbb C}/L_1,\dots,
E_{h_K}={\Bbb C}/L_{h_K}$ are pairwise non-isomorphic,
but their endomorphism ring is the same \cite{S}.


\end{document}